\documentclass[12pt]{amsart}
\usepackage{amscd, amsmath, amstext, amssymb, amsbsy}
\theoremstyle{plain}
\newtheorem{thm}{Theorem}[section]

\newtheorem{claim}[thm]{Claim}
\newtheorem{mydef}[thm]{Definition}

\newtheorem{cor}[thm]{Corollary}

\theoremstyle{definition}

\newcommand{\s}[1]{\mathcal{#1}}
\newcommand{\appa}{\subseteq}
\newcommand{\mbb}{\mathbb}
\newcommand{\im}{\mathrm{im} \, }

\newcommand{\spec}{\mathrm{spec} \, }

\newcommand{\lra}{\longrightarrow}

\newcommand{\sra}{\rightarrow}

\newcommand{\sotlim}{\text{SOT}-\lim}
\newcommand{\wotlim}{\text{WOT}-\lim}
\begin{document}
\bibliographystyle{amsplain}

\title{On free resolutions in multivariable 
operator theory} 
\author{Devin~C.V. Greene}
\address{Department of Mathematics\\ 
University of Nebraska\\ 
Lincoln, NE 68588} 
\email{dgreene@math.unl.edu} 


\begin{abstract}
Let $\s{A}_d$ be the complex polynomial ring in $d$ variables.  A 
contractive $\s{A}_d$-module is Hilbert space $\s{H}$ equipped with 
an $\s{A}_d$ action such that for any $\xi_1, \xi_2, \ldots, 
\xi_d \in \s{H}$, 

\begin{equation*}
\|z_1 \xi_1 + z_2 \xi + \cdots + z_d \xi_d \|^2 \leq \|\xi_1\|^2 
+ \|\xi_2\|^2 + \cdots + \|\xi_d\|^2.
\end{equation*}

\noindent Such objects have been shown to be useful for modeling 
$d$-tuples of 
mutually 
commuting operators acting on a Hilbert space.  An example is given 
by $H_d^2$, which is the Hilbert space of holomorphic functions on 
$B_d = \{z \in \mbb{C}^d : |z| \leq 1\}$  generated by the 
reproducing kernel $k(z,w) = \frac{1}{1-\left< z, w \right>}, \, z,w 
\in B_d$.  For a natural subcatagory of contractive 
$\s{A}_d$-modules, whose members are said to be {\it pure}, the 
module $H_d^2$ plays the role of the free module of rank one.  In 
fact, given any pure contractive $\s{A}_d$-module $\s{H}$, there is 
a ``free 
resolution'', i.e. an exact sequence of the following form:

\begin{equation*}
\begin{CD}
\cdots @>\Phi_2>> H_d^2 \otimes \s{C}_1 @>\Phi_1>> H_d^2 \otimes 
\s{C}_0 @>\Phi_0>> \s{H} @>>> 0.
\end{CD}
\end{equation*}

\noindent For $i \geq 1$, the map $\Phi_i$ can be viewed as 
a $\s{B}(\s{C}_i, \s{C}_{i-1})$-valued weakly holomorphic function 
on 
$B_d$.  ($\s{B}(\s{C}_i, \s{C}_{i-1})$ is the set of bounded linear 
operators from $\s{C}_i$ into $\s{C}_{i-1}$.)  ``Localizing'' the 
free resolution at a point $\lambda \in B_d$, one obtains a 
``localized complex'':

\begin{equation*}
\begin{CD}
\cdots @>\Phi_3(\lambda)>> \s{C}_2 @>\Phi_2(\lambda)>> \s{C}_1 
@>\Phi_1(\lambda)>> 
\s{C}_0.
\end{CD}
\end{equation*}

\noindent We shall show that the homology of this complex is 
isomorphic to the homology of the Koszul complex of the $d$-tuple 
$(\varphi^1, \varphi^2, \ldots, \varphi^d)$, where $\varphi^i$ is 
the $i$th coordinate function of a M\"obius transform on $B_d$ such 
that $\varphi(\lambda) = 0$.  

We shall also show that the set of M\"obius transforms on $B_d$ 
gives rise to a class of unitary operators on the Hilbert 
space $H_d^2$.  Explicitly, if $\varphi$ is a M\"obius transform on 
$B_d$, then the map $\xi \mapsto 
\frac{\sqrt{1-|\lambda|^2}}{1-\left< \cdot , \lambda \right>}(\xi 
\circ 
\varphi)$, where $\lambda = \varphi^{-1}(0)$, is a unitary 
operator on $H_d^2$.  We shall show that the set of such 
operators acts ``ergodically'' on $H_d^2$, in the sense that no 
non-trivial invariant subspace of $H_d^2$ is preserved by every 
operator of this form.  
\end{abstract}

\maketitle

\section{Introduction} \label{1}

Let $\s{A}_d = \mbb{C}[z_1, z_2, \ldots, z_d]$ be the polynomial ring 
in $d$ variables.  We define a
{\it contractive Hilbert $\s{A}_d$-module} to be a module $\s{H}$
over $\s{A}_d$ that is a Hilbert space and that has the additional
property that for all $\xi_1, \xi_2, \ldots, \xi_d \in \s{H}$, 

\begin{equation*} 
\|\sum_{k=1}^d \, z_k \xi_k\|^2 \leq \sum_{k=1}^d \, \|\xi_k\|^2.
\end{equation*}

\noindent Obviously any closed submodule $\s{K}$ of a contractive 
$\s{A}$-module $\s{H}$ is a contractive $\s{A}$-module.  Consider the 
Banach space quotient $\s{H}/\s{K}$.  One can identify this space 
with $\s{H} \ominus \s{K}$, and define a contractive $\s{A}_d$-module 
structure on it by compressing polynomials by $P_{\s{H} \ominus 
\s{K}}$.  

The notion of a Hilbert module was used by Arveson in 
\cite{arvcurvinv} to represent commuting $d$-contractions of 
operators.  Indeed, the actions of $z_1, z_2, \ldots, z_d$ on 
$\s{H}$ correspond to a mutually commuting $d$-tuple of 
linear operators 
$(T_1, T_2, \ldots, T_d)$ by defining $z_i \xi = T_ \xi$ for all 
$\xi \in \s{H}$.  The contractive condition on the module $\s{H}$ 
is equivalent to saying that the row operator $(T_1 \; T_2 \; \cdots 
\; T_d)$ is contractive when seen as a map from $\overbrace{\s{H} 
\oplus \s{H} \oplus \cdots \oplus \s{H}}^{d \text{~times~}} \lra 
\s{H}$.  Conversely, given a $d$-tuples of linear operators $(T_1, 
T_2, \ldots, T_d)$ acting on a Hilbert space $\s{H}$ where the 
components mutually commute and the row operator $(T_1 \; T_2 \; 
\cdots \; T_d)$ is contractive, one can define a contractive 
$\s{A}_d$-module structure on $\s{H}$ by defining $z_i \xi = T_i 
\xi$ for each $\xi \in \s{H}$.  Given a contractive 
$\s{A}_d$-module $\s{H}$, we call the $d$-tuple $(T_1, T_2, \ldots, 
T_d)$ the {\it associated $d$-tuple of $\s{H}$}.  

An example of a contractive $\s{A}_d$-module, and one which plays an 
important role in the 
theory, is the following:  Let $H_d^2$ be the Hilbert space of 
holomorphic functions on the unit ball $B_d$ of $\mbb{C}^d$ derived 
from 
the following reproducing kernel:

\begin{equation} \label{3}
k(z,\lambda) = \frac{1}{1-\left<z,\lambda\right>}, \quad z, \lambda 
\in B_d.
\end{equation}

\noindent The function $k_\lambda : z \mapsto k(z,\lambda)$ is in 
$H_d^2$ and in fact for any $\xi \in H_d^2$, 

\begin{equation*} 
\xi(\lambda) = \left< \xi, k_\lambda \right>.
\end{equation*}

\noindent A contractive $\s{A}_d$-module structure is defined on 
$H_d^2$ as follows:  For any $\xi \in H_d^2$, $p(z_1, z_2, \ldots, 
z_d) \in \s{A}_d$, 
$p(z_1, z_2, \ldots, z_d) \xi$ is simply the function $p(z_1, z_2, 
\ldots, z_d)$ multiplied by $\xi$.  The 
Hilbert space direct sum of $n$ copies of $H_d^2$ is also a 
contractive 
$\s{A}_d$-module and can be expressed by $H_d^2 \otimes \s{C}$, where 
$\s{C}$ is a Hilbert space of dimension $n$.  Note that we put no 
restrictions at this point on the cardinality of $n$.  For reasons 
that will soon become apparent, we call $H_d^2 \otimes \s{C}$ the 
{\it free contractive $\s{A}$-module of rank $\dim \s{C}$}.  

Another simple example of a contractive $\s{A}_d$-module arises as 
follows:  Let $C(\partial B_d)$ be the $C^*$-algebra of continuous 
functions on the unit sphere in $\mbb{C}^d$.  Let $\pi: 
C(\partial B_d) \lra \s{B}(\s{K})$ be a *-representation.  Then 
polynomials act on $\s{K}$ in the natural way, and $\s{K}$ becomes a 
contractive 
$\s{A}_d$-module.  Such modules are called {\it 
spherical modules}.  

Free and spherical contractive $\s{A}_d$-modules play the role of 
universal objects in the category of $\s{A}_d$-modules.  The 
following result 
is due to W. Arveson (\cite{arvmvot}).

\begin{thm} \label{4}
Let $\s{H}$ be a contractive $\s{A}_d$-module.  There exists a free 
module $\s{F}$, a spherical module $\s{S}$, and a module homomorphism 
$U: 
\s{F}
\oplus \s{S} \lra \s{H}$ such that $U U^* = 1$, i.e. $U$ is a 
coisometry.  
\end{thm}

\noindent A uniqueness condition also applies.  First we state the 
following definition: 

\begin{mydef} 
Let $\s{H}$ be a contractive $\s{A}_d$-module, let $\s{F}$ be a free 
module, let $\s{S}$ be a spherical module, and let $U:\s{F} \oplus 
\s{S} \lra 
\s{H}$ be a coisometry.  The triple $(\s{F},\s{S},U)$ is said to be a 
{\it 
minimal dilation} of $\s{H}$ if the closed submodule of $\s{F} \oplus 
S$ 
generated 
by $U^* \s{H}$ is $\s{F} \oplus \s{S}$.  
\end{mydef}

\noindent The following theorem, a proof of which can be found in 
\cite{arvcurvinv}, 
states that any two minimal dilations are naturally isomorphic.  

\begin{thm} \label{5.5}
Every contractive $\s{A}_d$-module $\s{H}$ has a minimal 
dilation.  Furthermore, if $(\s{F},\s{S},U)$ and $(\s{F}, 
\s{S'},U')$ are minimal dilations, 
then there is a 
unitary 
module isomorphism $V: \s{F} \oplus \s{S} \lra \s{F}' \oplus \s{S}'$ 
such 
that $U = 
U'V$.  Furthermore, $V$ has the form $V = (1_{H_d^2} \otimes W) 
\oplus W'$.  
\end{thm}

In this paper, we will be concerned almost entirely with {\it pure} 
contractive $\s{A}_d$-modules.  

\begin{mydef} \label{6}
A contractive $\s{A}_d$-module $\s{H}$ is {\it pure} if the spherical 
part of any minimal dilation is trivial.  In other words, $\s{H}$ is 
unitarily isomorphic to a quotient of a free module.  
\end{mydef}

\noindent There is an equivalent formulation of this definition
expressed in terms of contractive $d$-tuples of commuting operators.  
This equivalence follows from the work of Arveson in \cite{arvmvot}.

\begin{thm} \label{6.5}
Let $\s{H}$ be a contractive $\s{A}_d$-module, and let $T_1, T_2, 
\cdots, T_d$ be the linear operators on $\s{H}$ corresponding to the 
actions of $z_1, z_2, \ldots, z_d$.  Then $\s{H}$ is pure iff

\begin{equation*} 
\sotlim_{n \sra \infty} \, \sum_{i_1, i_2, \ldots, i_n}^d \, T_{i_1} 
T_{i_2} \cdots T_{i_n} T_{i_n}^* \cdots T_{i_2}^* T_{i_1}^* = 0.
\end{equation*}

\end{thm}

Let $\s{C}$ be a Hilbert space.  The free module $H_d^2 \otimes
\s{C}$ can be viewed as a space of (weakly) holomorphic
$\s{C}$-valued functions defined on $B_d$.  Indeed, if $\xi \in H_d^2
\otimes \s{C}$ and $\lambda \in B_d$, then $\xi(\lambda)$ is defined 
to be the unique element of $\s{C}$ such that for all $\eta \in 
\s{C}$, 

\begin{equation*} 
\left< \xi, k_\lambda \otimes \eta \right> = \left< \xi(\lambda), 
\eta \right>_{\s{C}}.
\end{equation*}

This identification of a free module with a space of vector valued 
holomorphic functions gives us a useful way of perceiving module 
homomorphisms between free modules.  Indeed, let $\s{D}$ and $\s{C}$ 
be Hilbert spaces, and let $\Phi: H_d^2 \otimes \s{D} \lra H_d^2 
\otimes \s{C}$ be a module homomorphism.  For each $\lambda \in 
B_d$, 
define the bounded linear operator $\Phi(\lambda): \s{D} \lra \s{C}$ 
by $\Phi(\lambda) \eta = \Phi(1 \otimes \eta)(\lambda)$.  Since 
$\Phi$ 
is a homomorphism, it follows that for all $\xi \in H_d^2 \otimes 
\s{D}$ and $\lambda \in B_d$, $(\Phi \xi)(\lambda) = \Phi(\lambda) 
\xi(\lambda)$.  Hence $\Phi$ is given by pointwise multiplication by 
the 
operator valued function $\Phi(\lambda)$.  From the definition, it is 
not difficult to show that the adjoint of $\Phi$ has the following 
property:  Let $\lambda \in B_d$, and let $\eta \in \s{C}$.  Then

\begin{equation} \label{8}
\Phi^*(k_\lambda \otimes \eta) = k_\lambda \otimes \Phi(\lambda)^* 
\eta.  
\end{equation}

\noindent As a consequence of \eqref{8}, we have the following 
boundedness condition on $\Phi(\lambda)$.

\begin{equation*} 
\|\Phi(\lambda)\| \leq \|\Phi\|,  \;  \forall \lambda \in B_d.
\end{equation*}

\noindent Thus the function $\lambda \in B_d \mapsto \Phi(\lambda)$
is a bounded $\s{B}(\s{D}, \s{C})$-valued (weakly) holomorphic
function on $B_d$.

Let $\s{H}$ be any pure contractive $\s{A}_d$-module, and let $\s{M}$
be a closed submodule.  Using Theorem~\ref{6.5}, it is not hard to
show that $\s{M}$ is pure as a contractive $\s{A}_d$-module.  Using
this, we can perform the following construction:  Let $\s{H}$ be a
pure contractive $\s{A}_d$-module.  Then by Theorem~\ref{4}, there
exists a free module $\s{F}_1$ and a coisometric module homomorphism
$\Phi_1: \s{F}_1 \lra \s{H}$.  Clearly $\ker \Phi_1$ is a closed
submodule of $\s{F}_1$, and by our above remarks, there exists a free
module $\s{F}_2$ and a coisometric module homomorphism $P: \s{F}_2
\lra \ker \Phi_1$.  Extending the codomain of $P$, we obtain a
partially isometric module homomorphism $\Phi_2: \s{F}_2 \lra 
\s{F}_1$ with range space $\ker \Phi_1$.  Continuing in this fashion, 
we obtain the following exact sequence:

\begin{equation} \label{10}
\begin{CD}
\cdots @>\Phi_2>> \s{F}_1 @>\Phi_1>> \s{F}_0 @>\Phi_0>> \s{H} @>>> 0,
\end{CD}
\end{equation}

\noindent where $\s{F}_i$ is free for each $i$.  We call this a 
free resolution of $\s{H}$, and it is 
analogous to the free resolution in the theory of finitely generated 
modules over Noetherian rings.  

Each $\s{F}_i$ appearing in \eqref{10} has the form $H_d^2 
\otimes 
\s{C}_i$ for some Hilbert space $\s{C}_i$.  If we localize to a point 
$\lambda \in B_d$, we obtain the following sequence of linear maps:

\begin{equation} \label{11}
\begin{CD}
\cdots @>\Phi_3(\lambda)>> \s{C}_2 @>\Phi_2(\lambda)>> \s{C}_1 
@>\Phi_1(\lambda)>> \s{C}_0
\end{CD}
\end{equation}

\noindent Since \eqref{10} is exact, it follows that 
\eqref{11} is a complex for each $\lambda$.  

The main result of Section~\ref{15} is that the homology of
\eqref{11} at $0$ is closely connected to the spectral properties of
the $d$-tuple $(T_1, T_2, \ldots, T_d)$ associated with $\s{H}$.  
By ``spectral properties'' we
are refering to the spectrum of a tuple of operators in the sense of 
Taylor (cf. \cite{taylor1}).
The following discussion summarizes this more precisely:  For each 
$k \in \mbb{Z}$, let $\bigwedge^k \, \mbb{C}^d$ be the $k$-fold 
alternating product of $\mbb{C}^d$, taken to be the trivial vector 
space when $k < 0$.  For each $k$, let $E_k = \s{H} \otimes 
\bigwedge^k \, \mbb{C}^d$.  Fix an orthonormal basis $\{e_1, e_2, 
\ldots, e_d\}$ for $\mbb{C}^d$, and for each $k$, define the linear 
map $\partial_k : E_k \lra E_{k-1}$ as follows:

\begin{equation*} 
\partial_k (\xi \otimes e_{i_1} \wedge \cdots \wedge e_{i_k}) = 
\sum_{j=1}^k \, (-1)^{j+1} z_{i_j} \xi \otimes e_{i_1} \wedge \cdots 
\wedge 
\hat{e_{i_j}} \wedge \cdots \wedge e_{i_k}.
\end{equation*}

\noindent A straightforward calculation shows that the following 
sequence is a complex.

\begin{equation} \label{13}
\begin{CD}
0 @>>> E_d @>\partial_d>> E_{d-1} @>\partial_{d-1}>> 
\cdots @>\partial_1>> E_0 @>>> 0 
\end{CD}
\end{equation}

\noindent What we will show in Section~\ref{15} is that 

\begin{equation*} 
\frac{\ker \partial_k}{\im \partial_{k+1}} \cong \frac{\ker 
\Phi_k(0)}{\im \Phi_{k+1}(0)}.
\end{equation*}

\noindent Hence the homology of \eqref{13} and that of 
\eqref{11} when $\lambda = 0$ are identical.  

In Section~\ref{60}, we introduce the notion of a M\"obius transform 
of 
a pure contractive $\s{A}_d$-module.  We recall that a M\"obius 
transform $\varphi$ on $B_1$ is a bijective holomorphism from $B_d$ 
onto itself.  
If $\s{H}$ is a pure contractive $\s{A}_d$-module, and 
$(T_1, T_2, \ldots, T_d)$ is the associated $d$-tuple,  then we can 
define the
M\"obius transform $(\s{H})_\varphi$ of $\s{H}$ to be the pure
contractive $\s{A}_d$-module where the underlying Hilbert space of
$\s{H}_\lambda$ is the underlying Hilbert space of $\s{H}$, but $z_1, 
z_2, \ldots,
z_d$ act as $\varphi^1, \varphi^2, \ldots, \varphi^d$, respectively, 
where 
$\varphi^i$ is the $i$th coordinate of $\varphi$.  We will show
that there exists a unitary module isomorphism $U_\varphi: H_d^2 \lra
(H_d^2)_\lambda$, which carries the space of all functions in $H_d^2$
that vanish at $0$ to the space of all functions that vanish at
$\lambda$.  We will also show that the $U_\lambda$'s act 
``ergodically'' on $H_d^2$ in the sense that no proper nontrivial 
closed submodule of 
$H_d^2$ is invariant under $U_\lambda$ for all $\lambda \in B_d$.  

Finally, in Section~\ref{85}, we use the results on Section~\ref{60} 
to 
describe the
homology of \eqref{11}
for an arbitrary $\lambda \in B_d$.  In particular, Theorem~\ref{87} 
states that the homology of the Koszul complex of $(\s{H})_\varphi$
is equivalent to the localized
complex \eqref{11}.  Hence \eqref{11} contains spectral
information (in the sense of Taylor) of the M\"obius transformation
of $\s{H}$.

\section{The Koszul complex and free resolution of a contractive 
$\s{A}_d$-module} \label{15}

In his seminal paper (\cite{taylor1}), J. Taylor defined the notion of 
invertibility of $d$-tuples of commuting operators acting on a Banach 
space $\s{B}$.  We will briefly summarize Taylor's construction:

Let $\s{B}$ be a Banach space, and let $(a_1, a_2, \ldots a_d)$ be a 
$d$-tuple of mutually commuting bounded operators on $\s{B}$.  For $k 
\in \mbb{Z}$, let $\bigwedge^k \, \mbb{C}^d$ be the $k$-fold wedge 
product of $\mbb{C}^d$, where this is taken to be the trivial vector 
space when $k < 0$.  Let $E_k = \s{B} \otimes \bigwedge^k \, 
\mbb{C}^d$.  Note that $E_k$ can be viewed as an $\left( 
\begin{array}{c} 
n \\
k \end{array} \right)$-fold direct sum of $E_k$'s, hence it is a 
Banach 
space itself.  Fix an orthonormal basis $\{e_1, e_2, \ldots e_d\}$ 
for $\mbb{C}^d$.  For each $k$, define $\partial_k : E_k \lra 
E_{k-1}$ by 

\begin{equation*} 
\partial_k \xi \otimes e_{i_1} \wedge \cdots \wedge e_{i_k} = 
\sum_{j=1}^d \, (-1)^{j+1} a_i \xi \otimes e_{i_1} \wedge \cdots 
\hat{e_{i_j}} 
\cdots e_{i_k},
\end{equation*}

\noindent A simple computation shows that $\partial_{k-1} 
\partial_k = 0$.  Hence we have the following complex of Banach 
spaces.  

\begin{equation} \label{17}
\begin{CD}
\cdots @>>> 0 @>>> E_d @>\partial_d>> E_{d-1} @>\partial_{d-1}>> 
\cdots @>\partial_1>> E_0 @>\partial_0>> 0.
\end{CD}
\end{equation}

\noindent The following is Taylor's definition of invertibility:

\begin{mydef} \label{17.5}
A $d$-tuple of mutually commuting operators acting on a common Banach 
space $\s{B}$ is said to be invertible if the sequence in \eqref{17} 
is exact.  
\end{mydef}

\noindent One can immediately provide a definition of the spectrum 
of $(a_1, a_2,$ $\ldots, a_d)$:

\begin{mydef} \label{18}
Let $a = (a_1, a_2, \ldots, a_d)$ be a $d$-tuple of mutually 
commuting operators on a Banach space $\s{B}$.  Then $\spec (a)$ is 
defined to be the set $\lambda = (\lambda_1, \lambda_2, \ldots, 
\lambda_d) \in \mbb{C}^d$ where $(a_1 - \lambda_1, a_2 - \lambda_2, 
\ldots, a_d -\lambda_d)$ is not invertible.  
\end{mydef}

We remark that in the case where $d=1$, both Definition~\ref{17.5} 
and Definition~\ref{18} reduce to the usual definitions of 
invertibility and spectrum for single operators.  

If we restrict ourselves to the case where $\s{B}$ is a 
Hilbert
space $\s{H}$, and the $d$-tuple $(T_1, T_2, \ldots, T_d)$ consists
of mutually commuting bounded operators on $\s{H}$, then we can 
express
homological features of \eqref{17} in terms of a single self-adjoint 
operator.  Our presentation will follow Arveson (\cite{arvdirac}),
but this idea appeared earlier in the work of such authors as Curto 
(\cite{curto1})
and Vasilescu (\cite{vasilescu1}).

Let $\s{H}$ and $(T_1, T_2, \ldots, T_d)$ be as above, and let$\{e_1, 
e_2, \ldots, e_d\}$ be our fixed orthonormal basis for $\mbb{C}^d$.  
The space $\bigwedge^k \, \mbb{C}^d$ has a natural inner product 
defined by 

\begin{equation*} 
\left< z_1 \wedge z_2 \wedge \cdots \wedge z_k, w_1 \wedge w_2 \wedge 
\cdots \wedge w_k \right> = \det (\left< z_i, w_j \right>)_{ij}, 
\; z_i, w_i \in \mbb{C}^d.
\end{equation*}

\noindent Hence the spaces $E_k$ are tensor products of Hilbert 
spaces.  Let $\bigwedge \, \mbb{C}^d = \bigoplus_{k \in \mbb{Z}} \, 
\bigwedge^k \, \mbb{C}^d$.  Define $E$ to be the direct sum of the 
$E_k$'s, i.e.

\begin{equation*} 
E = \bigoplus_{k \in \mbb{Z}} \, E_k = \s{H} \otimes \bigwedge \, 
\mbb{C}^d.
\end{equation*}

\noindent Let $c_1, c_2, \ldots, c_d$ be the operators on $\bigwedge 
\, \mbb{C}^d$ defined by 

\begin{equation*} 
c_i (z_1 \wedge z_2 \wedge \cdots \wedge z_k) = e_i \wedge z_1 \wedge 
\cdots \wedge z_k, \; z_i \in \mbb{C}^d.  
\end{equation*}

\noindent We then define the linear operator $\partial: E \lra E$ to 
be the sum 

\begin{equation*} 
T_1 \otimes c_1^* + T_2 \otimes c_2^* + \cdots + T_d \otimes 
c_d^*.
\end{equation*}

\noindent One then checks that the restriction of $\partial$ to $E_k$ 
is $\partial_k$.  One can now prove the following theorem (cf. 
\cite{arvdirac}):

\begin{thm} \label{23}
Let $(T_1, T_2, \ldots, T_d)$ be a $d$-tuple of mutually commuting 
operators on a common Hilbert space $\s{H}$.  If $\partial$ is the 
corresponding boundary operator, then $(T_1, T_2, \ldots, T_d)$ is 
invertible iff $\partial + \partial^*$ is invertible.  
\end{thm}

The operator $\partial + \partial^*$ is called the Dirac operator 
corresponding to $(T_1, T_2, \ldots, T_d)$.  This 
idea of expressing the invertibility of a $d$-tuple in terms of 
a single operator suggests the following definition:

\begin{mydef} 
A $d$-tuple of operators as in Theorem~\ref{23} is said to be 
Fredholm if the corresponding Dirac operator $\partial + 
\partial^*$ is Fredholm.  
\end{mydef}

\noindent We note that in the case where $d=1$, this definition 
corresponds to the usual definition of Fredholmness via an easy 
application of Atkinsons' equivalences (see \cite{conway90course}, 
for 
example).  

The Fredholmness of a Dirac operator has an importance consequence 
with respect to the homology of the Koszul complex.  This is 
expressed in the following theorem, which follows from the definition 
of $\partial$ by a 
straightforward argument

\begin{thm} \label{25} 
A $d$-tuple $(T_1, T_2, \ldots, T_d)$ is Fredholm iff the
homology spaces $\ker \partial_k/\im \partial_{k+1}$ are finite
dimensional for all $k \in \mbb{Z}$.  
\end{thm}

\noindent Theorem~\ref{25} allows us to generalize the notion of 
index to Fredholm $d$-tuples.  Indeed, the index of a Fredholm 
$d$-tuple $(T_1, T_2, \ldots, T_d)$ is defined to be the alternating 
sum of the dimensions of the homology spaces of the Koszul complex, 
i.e.

\begin{equation*} 
\text{index} \, (T_1, T_2, \ldots, T_d) = \sum_{k \in \mbb{Z}} \, 
(-1)^{k+1} \, \frac{\ker \partial_k}{\im \partial_{k+1}}.
\end{equation*}

\noindent Again, we note that this definition reduces to the usual 
definition of index when one takes $d=1$.  

We now show that the $d$-tuple associated with the free Hilbert
module $H_d^2$ is Fredholm and we compute its homology.

\begin{thm} 
The $d$-tuple $(S_1, S_2, \ldots, S_d)$ associated with $H^2_d$ is 
Fredholm.  Furthermore, the extended sequence of maps

\begin{equation*} 
\begin{CD}
\cdots @>\partial_2>> E_1 @>\partial_1>> E_0 @>\partial_{1/2}>> 
\mbb{C} @>>> 0,
\end{CD}
\end{equation*}

\noindent where $\partial_{1/2}$ is defined to be the evalation map 
$\xi
\mapsto \xi(0)$, is exact.  
\end{thm}

\noindent A proof of the fact the $(S_1, S_2, \ldots, S_d)$ is 
Fredholm can be found in \cite{arvdirac}.  We will rely on this fact 
to 
prove the remainder of the theorem.  

\begin{proof} 
Let $k$ be an integer no less than $1$.  By Theorem~\ref{25}, the 
space 
$\im \partial_{k+1}$ has finite codimension in $\ker \partial_k$.  It 
is a standard fact in operator theory that if the image of a bounded 
operator has finite codimension in a larger closed subspace, then the 
image is closed.  Hence $\im \partial_{k+1}$ is a closed subspace of 
finite codimension in $\ker \partial_k$.  We now show that $\im 
\partial_{k+1} = \ker \partial_k$.  To this end, let $\xi \in \ker 
\partial_k$.  Then $\xi$ can be written as follows:

\begin{equation*} 
\sum_{1 \leq i_1 < i_2 < \ldots <  i_k  \leq d} 
\xi_{i_1, i_2, 
\ldots, i_k} \otimes e_{i_1} \wedge e_{i_2} \wedge \cdots \wedge 
e_{i_k}
\end{equation*}

\noindent where $\xi_{i_1, i_2, \ldots i_k} \in H_d^2$.  Supposing 
for the moment that these $\xi_{i_1, i_2, \ldots, i_k}$ are 
homogeneous 
polynomials all of the same degree $N$, then $\partial_k \xi$ is in 
a similar form but with common degree $N + 1$. It follows in the 
general case that if $\xi_{i_1, i_2, \ldots, i_k}^n$ is the $n$th 
degree 
homogeneous component of $\xi_{i_1, i_2, \ldots, i_k}$, then 

\begin{equation*} 
\xi^n := \sum_{1 \leq i_1 < i_2 < \cdots < i_k \leq d} \, \xi_{i_1, 
i_2, 
\ldots, i_k}^n \otimes e_{i_1} \wedge e_{i_2} \wedge \cdots \wedge 
e_{i_k} \in \ker \partial_k.
\end{equation*}

\noindent By Corollary~17.5 in \cite{eisenbud1}, $\xi^n \in \im 
\partial_{k+1}$ for each 
$n$, hence $\xi = \sum_{n=0}^\infty \, \xi_n \in \im \partial_{k+1}$.  

We proceed to the case where $k=0$.  By Theorem~\ref{68},
the following row operator is a partial isometry with image $\{\xi 
\in H_d^2 : \xi(0) = 0\}$.  

\begin{equation} \label{31}
(S_1 \; S_2 \; \cdots \; S_d): \overbrace{H_d^2 
\oplus \cdots \oplus H_d^2}^{d \text{~times~}} \lra H_d^2.
\end{equation}

\noindent A generic element $\zeta$ of $E_1$ has the following form 

\begin{equation*} 
\sum_{k=1}^d \, \xi_k \otimes e_k, \; \xi_k \in H_d^2.
\end{equation*}

\noindent Hence $\partial_1 \xi = \sum_{k=1}^d \, S_k \xi_k$.  It 
follows from \eqref{31} and the statement preceding it that $\im 
\partial_1 = \ker \partial_{1/2}$.  

The surjectivity of $\partial_{1/2}$ is clear.
\end{proof}

\begin{cor} 
Let

\begin{equation*} 
\begin{CD}
\cdots @>\partial_2>> E_1 @>\partial_1>> E_0 @>>> 0
\end{CD}
\end{equation*}

\noindent be the Koszul complex of $(S_1, S_2, \cdots, S_d)$.  Then 
its homology is as follows:

\begin{equation*} 
\frac{\ker \partial_k}{\im \partial_{k+1}} = 0, \; k
\geq 
1 
\end{equation*}

\begin{equation*}
\frac{\ker \partial_0}{\im \partial_{1}} \cong 
\mbb{C}.
\end{equation*}

\end{cor}

\begin{cor} 
Let $(S'_1, S'_2, \ldots, S'_d)$ be the $d$-tuple associated with a 
free module $\s{F} = H_d^2 \otimes \s{C}$ with the following extended 
sequence of maps:

\begin{equation*} 
\begin{CD}
\cdots @>\partial_2>> E_1 @>\partial_1>> E_0 @>\partial_{1/2}>> \s{C} 
@>>> 0.
\end{CD}
\end{equation*}

\noindent Then for $k \geq 1$, $\ker \partial_k = \im 
\partial_{k+1}$
and $\ker \partial_0/\im \partial_{1} \cong \ker 
\partial_0 \cap (\im \partial_{1})^\perp \cong \s{C}$.  In other 
words $\im \partial_{k+1}$ is closed in $\ker \partial_k$, and in the 
case where $k=1$, the codimension of $\im \partial_{k+1}$ in $\ker 
\partial_k$ is $\dim \s{C}$.  

\end{cor}

Naturally, our entire discussion on Koszul complexes of $d$-tuples of
operators can be rephrased in terms of $\s{A}_d$-modules.  Indeed,
one simply takes $\s{H}$ to be the module defined by $z_i \xi = T_i
\xi$ for all $\xi \in \s{H}$.  The spaces $E_k$ are defined 
analogously, with $\partial_k : E_k \lra E_{k-1}$ reexpressed as 

\begin{equation*} 
\partial_k (\xi \otimes e_{i_1} \wedge e_{i_2} \wedge \cdots \wedge 
e_{i_k}) = \sum_{j=1}^k \, (-1)^{j+1} z_{i_j} \xi \otimes e_{i_1} 
\wedge \cdots 
\wedge \hat{e_{i_j}} \wedge \cdots \wedge e_{i_k}, 
\end{equation*} 

\noindent Since $\s{A}_d$ is a commutative algebra, the maps 
$\partial_k$ are module homomorphisms.  Hence the sequence \eqref{17} 
can be viewed as a complex of $\s{A}_d$-modules.  

As summarized in Section~\ref{1},
starting with a pure $\s{A}_d$-module $\s{H}$, by means of dilation 
theory one may construct a free resolution of $\s{H}$:

\begin{equation} \label{38}
\begin{CD}
\cdots @>\Phi_2>> \s{F}_1 @>\Phi_1>> \s{F}_0 @>\Phi_0>> \s{H} @>>> 0, 
\end{CD}
\end{equation}

\noindent which is an exact sequence where the $\Phi_i$'s are partial
isometries and the $\s{F}_i$'s are free modules of the form $H_d^2
\otimes \s{C}_i$.  One then ``localizes'' \eqref{38} to a point
$\lambda \in B_d$ to obtain a complex of vector spaces:

\begin{equation*} 
\begin{CD}
\cdots @>\Phi_3 (\lambda)>> \s{C}_2 @>\Phi_2 (\lambda)>> \s{C}_1 
@>\Phi_1 (\lambda)>> \s{C}_0.
\end{CD}
\end{equation*}

\noindent The main result of this section is the following:

\begin{thm} \label{39.25}
Let $\s{H}$ be a pure contractive $\s{A}_d$-module, and let 

\begin{equation} \label{39.5}
\begin{CD}
\cdots @>\Phi_2>> \s{F}_1 @>\Phi_1>> \s{F}_0 @>\Phi_0>> \s{H} @>>> 0.
\end{CD}
\end{equation}

\noindent a free resolution of $\s{H}$.  Localize at $0$ to obtain the 
following complex:

\begin{equation*} 
\begin{CD}
\cdots @>\Phi_3(0)>> \s{C}_2 @>\Phi_2(0)>> \s{C}_1 @>\Phi_1(0)>> 
\s{C}_0.
\end{CD}
\end{equation*}

\noindent Let 

\begin{equation*} 
\begin{CD}
\cdots @>\partial_3>> E_2 @>\partial_2>> E_1 @>\partial_1>> E_0 @>>> 
0
\end{CD}
\end{equation*}

\noindent be the Koszul complex of $\s{H}$.  Then for all $k \geq 1$, 

\begin{equation*} 
\frac{\ker \partial_k}{\im \partial_{k+1}} \cong \frac{\ker 
\Phi_{k}(0)}{\im \Phi_{k+1}(0)}.
\end{equation*}

\end{thm}

\begin{proof}
In the course of this proof, we will use $\s{F}_i^k$ and $\s{H}^k$ to 
denote, respectively, $\s{F}_i \otimes \bigwedge^k \, \mbb{C}^d$ and 
$\s{H} \otimes \bigwedge^k \, \mbb{C}^d$.  Since for any $k$, 
$\bigwedge^k \, \mbb{C}^d$ is finite dimensional, tensoring the 
components of \eqref{39.5} by $\bigwedge^k \, \mbb{C}^d$ preserves 
exactness.  Hence the following sequence is exact:

\begin{equation*} 
\begin{CD}
\cdots @>\Phi_2 \otimes 1_{\bigwedge^k \, \mbb{C}^d}>> \s{F}_1^k
@>\Phi_1 \otimes 1_{\bigwedge^k \, \mbb{C}^d}>> \s{F}_0^k @>\Phi_0
\otimes 1_{\bigwedge^k \, \mbb{C}^d}>> \s{H}^k @>>> 0.
\end{CD}
\end{equation*}

\noindent For the sake of convenience, and since it will cause no
confusion in what follows, we will denote maps of the form $A \otimes 
\bigwedge^k \, \mbb{C}^d$ by $A$ for any linear operator $A$ on 
$H_d^2$.  Furthermore, unless otherwise stated, we shall 
use $\partial_k$ to denote {\it any} Koszul complex mapping $\s{M} 
\otimes \bigwedge^k \, \mbb{C}^d \lra \s{M} \otimes \bigwedge^{k-1} \, 
\mbb{C}^d$.  

We claim that the following diagram commutes, and, with the 
exception of the $\s{C}_*$ column, is exact on rows and columns.  

\begin{equation} \label{43.25}
\begin{CD}
& & \vdots & & \vdots & & \vdots\\	
& & @V\Phi_{i+3}VV @V\Phi_{i+3}VV @V\Phi_{i+3}(0)VV \\
\cdots @>\partial_2>> \s{F}_{i+2}^1 @>\partial_1>> \s{F}_{i+2}^0 
@>\partial_{1/2}>> \s{C}_{i+2} @>>> 0 \\
& & @V\Phi_{i+2}VV @V\Phi_{i+2}VV @V\Phi_{i+2}(0)VV \\
\cdots @>\partial_2>> \s{F}_{i+1}^1 @>\partial_1>> \s{F}_{i+1}^0
@>\partial_{1/2}>> \s{C}_{i+1} @>>> 0 \\
& & @V\Phi_{i+1}VV @V\Phi_{i+1}VV @V\Phi_{i+1}(0)VV \\
\cdots @>\partial_2>> \s{F}_{i}^1 @>\partial_1>> \s{F}_{i}^0
@>\partial_{1/2}>> \s{C}_{i} @>>> 0 \\
& & @V\Phi_{i}VV @V\Phi_{i}VV @V\Phi_{i}(0)VV \\
& & \vdots & & \vdots & & \vdots\\ 
\end{CD}
\end{equation}

\noindent First we check commutativity for the square 

\begin{equation} \label{43.5}
\begin{CD}
\s{F}_{i+1}^k @>\partial_k>> \s{F}_{i+1}^{k-1} \\
@V\Phi_{i+1}VV @V\Phi_{i+1}VV \\
\s{F}_i^k @>\partial_k>> \s{F}_i^{k-1},
\end{CD}
\end{equation}

\noindent where $i,k \geq 1$.  This amounts to showing that 
$\partial_k \Phi_{i+1} = 
\Phi_{i+1} \partial_k$.  Consider an element in $\s{F}_{i+1}^k$ 
of the 
form 

\begin{equation} \label{44}
\xi \otimes e_{i_1} \wedge e_{i_2} \wedge \cdots \wedge e_{i_k}, 
\end{equation}

\noindent where $\xi \in \s{F}_{i+1}, 1 \leq i_1 < i_2 < \cdots < 
i_k 
\leq d$.  Applying $\Phi_{i+1}$ and then $\partial_k$ to this gives 
us

\begin{equation} \label{45}
\sum_{j=1}^k \, (-1)^{j+1} z_{i_j} \Phi_{i+1} \xi \otimes e_{i_1} 
\wedge \cdots 
\wedge \hat{e_{i_j}} \wedge \cdots e_{i_k}.
\end{equation}

\noindent Since $\Phi_{i+1}$ is a module homomorphism, 
the $z_{i_j}$'s 
commute with $\Phi_{i+1}$, and the resulting expression 

\begin{equation*} 
\sum_{j=1}^k \, (-1)^{j+1} \Phi_{i+1} z_{i_j} \xi \otimes 
e_{i_1} \wedge \cdots
\wedge \hat{e_{i_j}} \wedge \cdots e_{i_k}
\end{equation*}

\noindent is the result of applying $\partial_k$ then $\Phi_{i+1}$ to 
\eqref{44}.  Hence \eqref{43.5} is a commuting square.  For the case 
of squares of the following form:

\begin{equation} \label{47}
\begin{CD}
\s{F}_{i+1}^0 @>\partial_{1/2}>> \s{C}_{i+1} \\
@V\Phi_{i+1}VV @V\Phi_iVV \\
\s{F}_i^0 @>\partial_{1/2}>> \s{C}_i,
\end{CD}
\end{equation}

\noindent we argue as follows.  Let $\xi \in \s{F}_{i+1}^0 = 
\s{F}_{i+1}$.  The result of 
evaluating at $0$ and then applying $\Phi_{i+1}(0)$ results in 

\begin{equation*} 
\Phi_{i+1}(0) \xi(0)
\end{equation*}

\noindent which is equivalent to $(\Phi_{i+1} \xi)(0)$, which is the 
result of 
applying $\Phi_{i+1}$ to $\xi$ and then evaluating at $0$.  Hence 
\eqref{47} is a commuting square.  

Before proceding we make the following observation:  For each $i \geq 
0$, we denote the quotient module $\s{F}_i/\Phi_{i+1}(\s{F}_{i+1})$ 
by 
$\s{H}_i$.  Note that by assumption $\s{H}$ is naturally isomorphic to 
$\s{H}_0$.  Let 

\begin{equation} \label{49}
\begin{CD}
\cdots @>\partial'_2>> \s{H}_i^1 @>\partial'_1>> \s{H}_i^0 @>>> 0 
\end{CD}
\end{equation}

\noindent be the Koszul complex for $\s{H}_i$.  Fixing this $i$, we 
let 

\begin{equation} \label{50} 
\begin{CD}
\cdots @>\partial''_2>> \s{H}_{i+1}^1 @>\partial''_1>> \s{H}_{i+1}^0 
@>>> 0
\end{CD}
\end{equation}

\noindent be the Koszul complex for $\s{H}_{i+1}$.  Note that since 
\eqref{43.25} commutes, the complexes \eqref{49} and \eqref{50} are 
induced by the maps of the Koszul complexes $\s{F}_{i+1}^* \lra 
\s{F}_i^*$ and $\s{F}_{i+2}^* \lra \s{F}_{i+1}^*$ respectively.  
This can be expressed by saying that the following two diagrams
commute:

\begin{equation*} \label{50.5}
\begin{CD} 
\cdots  @>\partial_2>> \s{F}_{i+1}^1 @>\partial_1>> \s{F}_{i+1}^0 \\
& & @V\Phi_{i+1}VV @V\Phi_{i+1}VV \\
\cdots @>\partial_2>> \s{F}_{i}^1 @>\partial_1>> \s{F}_{i}^0 \\
& & @VVV @VVV \\ 
\cdots @>\partial'_2>> \s{H}_i^1 @>\partial'_1>> \s{H}_i^0 \\
& & @VVV @VVV \\
& & 0 & & 0,
\end{CD}
\end{equation*}

\vspace{5mm}

\begin{equation*}
\begin{CD}
\cdots @>\partial_2>> \s{F}_{i+2}^1 @>\partial_1>> \s{F}_{i+2}^0\\
& & @V\Phi_{i+2}VV @V\Phi_{i+2}VV \\
\cdots @>\partial_2>> \s{F}_{i+1}^1 @>\partial_1>> \s{F}_{i+1}^0 \\
& & @VVV @VVV \\
\cdots @>\partial''_2>> \s{H}_{i+1}^1 @>\partial''_1>> \s{H}_{i+1}^0 
\\
& & @VVV @VVV \\ 
& & 0 & & 0
\end{CD}
\end{equation*}

\begin{equation*} 
\begin{CD}
\end{CD}
\end{equation*}

\noindent We will use this fact to establish the following two claims:  

\begin{claim} \label{51}
For $k \geq 2$, 

\begin{equation*} 
\frac{\ker \partial'_k}{\im \partial'_{k+1}} \cong \frac{\ker 
\partial''_{k-1}}{\im \partial''_k}.
\end{equation*}

\end{claim}

\begin{proof}
Consider the following portion of \eqref{43.25}, where $k \geq 2$:

\begin{equation} \label{53}
\begin{CD}
\s{F}_{i+2}^{k+1} @>\partial_{k+1}>> \s{F}_{i+2}^k @>\partial_k>> 
\s{F}_{i+2}^{k-1} @>\partial_{k-2}>> \s{F}_{i+2}^{k-2}\\
@V\Phi_{i+2}VV @V\Phi_{i+2}VV @V\Phi_{i+2}VV  @V\Phi_{i+2}VV\\
\s{F}_{i+1}^{k+1} @>\partial_{k+1}>> \s{F}_{i+1}^k @>\partial_k>>
\s{F}_{i+1}^{k-1} @>\partial_{k-1}>> \s{F}_{i+1}^{k-2} \\
@V\Phi_{i+1}VV @V\Phi_{i+1}VV @V\Phi_{i+1}VV @V\Phi_{i+1}VV \\
\s{F}_{i}^{k+1} @>\partial_{k+1}>> \s{F}_{i}^k @>\partial_k>>
\s{F}_{i}^{k-1} @>\partial_{k-1}>> \s{F}_i^{k-2}.
\end{CD}
\end{equation}

\noindent Let $\zeta \in \ker \partial'_k/\im \partial'_{k+1}$.  
Choose a representative $\zeta_0 \in \s{F}_i^k$ for $\zeta$.  By
assumption, $\partial_k \zeta_0 \in \Phi_{i+1}(\s{F}_{i+1}^{k-1})$.  
Choose $\eta_0 \in \s{F}_{i+1}^k$ such that $\Phi_{i+1} \eta_0 =
\partial_k \zeta_0$.  By exactness of the bottom row of \eqref{53}, 
$\partial_{k-1} \Phi_{i+1} \eta_0 = 0$, hence $\partial_{k-1} \eta_0 
\in \ker \Phi_{i+1}$.  Hence by the exactness of the last column in 
\eqref{53}, $\eta_0$ defines a homology class $\eta \in \ker 
\partial''_{k-1}/\im \partial''_k$.  We claim that 
this $\eta$ depends only on the choice of $\zeta$.  Indeed, due to 
exactness of rows in \eqref{53}, a 
different choice of $\eta_0$ corresponds to a perturbation by 
an element in the image of $\s{F}_{i+2}^{k-1}$ under $\Phi_{i+2}$, 
which results in the new $\eta_0$ being in the same homology class.  A 
different choice of $\zeta_0$ corresponds to a perturbation by an 
element of $\partial_{k+1} (\s{F}_i^{k+1})$ and an element of 
$\Phi_{i+1}(\s{F}_{i+1}^k)$.  The first of these is eliminated by the 
exactness of the bottom row of \eqref{53} and the way we defined 
$\eta$.  Due to the commutativity of the diagram, the second 
pertubation corresponds to a perturbation of $\eta_0$ by an element of 
$\partial_k (\s{F}_{i+1}^k)$, which yield the same homology class 
$\eta$.  

Conversely, if one begins with an element $\eta \in \ker
\partial''_{k-1}/\im \partial''_k$, for any representative $\eta_0 \in
\s{F}_{i+1}^{k-1}$, there is an element $\zeta_0 \in \s{F}_i^k$ such
that $\partial_k \zeta_0 = \Phi_{i+1} \eta_0$.  This implies that 
$\zeta_0$ corresponds to a homology class $\zeta \in \ker 
\partial'_k/\im \partial'_{k+1}$.  We claim that this $\zeta$ depends 
only on $\eta$.  Indeed, by the exactness of the bottom row, a 
different choice of $\zeta_0$ is a perturbation by an element in 
$\partial_{k+1} (\s{F}_i^{k+1})$, which corresponds to the same 
homology class $\zeta$.  A different choice of $\eta_0$ corresponds to 
a perturbation by an element in $\Phi_{i+2}(\s{F}_{i+2}^{k-1})$ or in 
$\partial_k (\s{F}_{i+1}^k)$.  In the first case, the perturbation is 
annihilated by the time that we get to $\s{F}_i^{k-1}$.  In the second 
case, the commutativity of the diagram implies that the resulting 
perturbation of $\zeta_0$ is an element in 
$\Phi_{i+1}(\s{F}_{i+1}^k)$ and hence the homology class $\zeta$ is 
the same.  

Clearly the above maps are inverses of each other, and hence we have 
the desired isomorphism.  

\end{proof}

\begin{claim} \label{54}
\begin{equation*} 
\frac{\ker \partial'_1}{\im \partial'_2} \cong \frac{\ker 
\Phi_{i+1}(0)}{\im \Phi_{i+2}(0)}.
\end{equation*}
\end{claim}

\begin{proof}
As in Claim~\ref{51}, we focus on a piece of \eqref{43.25}:

\begin{equation} \label{55}
\begin{CD}
\s{F}_{i+2}^2 @>\partial_2>> \s{F}_{i+2}^1 @>\partial_1>> 
\s{F}_{i+2}^0 @>\partial_{1/2}>> \s{C}_{i+2} @>>> 0 \\
@V\Phi_{i+2}VV @V\Phi_{i+2}VV @V\Phi_{i+2}VV @V\Phi_{i+2}(0)VV \\
\s{F}_{i+1}^2 @>\partial_2>> \s{F}_{i+1}^1 @>\partial_1>>
\s{F}_{i+1}^0 @>\partial_{1/2}>> \s{C}_{i+1} @>>> 0 \\ 
@V\Phi_{i+1}VV @V\Phi_{i+1}VV @V\Phi_{i+1}VV @V\Phi_{i+1}(0)VV \\ 
\s{F}_{i}^2 @>\partial_2>> \s{F}_{i}^1 @>\partial_1>>
\s{F}_{i}^0 @>\partial_{1/2}>> \s{C}_{i} @>>> 0.
\end{CD}
\end{equation} 

\noindent Let $\zeta$ be an element in $\ker \partial'_1/\im 
\partial'_{2}$,  
let $\zeta_0$ be a representative of $\zeta$ in $\s{F}_i^1$.  Then 
$\partial_1 \zeta_0 \in \Phi_{i+1}(\s{F}_{i+1}^0)$ by assumption.  
Choose $\eta \in \s{F}_{i+1}^0$ such that $\Phi_{i+1} \eta = 
\partial_1 \zeta_0$, and then let $x_0 = \eta(0) \in \s{C}_{i+1}$.  By 
the way in which $x_0$ was defined and by the commutativity of 
\eqref{55}, $x_0 \in \ker \Phi_{i+1}(0)$.  Now let $x$ be the 
homology class of $x_0$ in $\ker \Phi_{i+1}(0)/\im \Phi_{i+2}(0)$.  We 
claim that this $x$ depends only on $\zeta$.  First, by exactness, a 
different choice of $\eta$ corresponds to a perturbation by an element 
in $\Phi_{i+2}(\s{F}_{i+2}^0)$.  By the commutativity of \eqref{55}, 
this corresponds to a perturbation by $x_0$ by an element in 
$\Phi_{i+2}(0)(\s{C}_{i+2})$, which leaves $x$ unaltered.  By the 
exactness of rows, any perturbation of $\zeta_0$ by an element of 
$\partial_2 (\s{F}_i^2)$ or $\Phi_{i+1} (\s{F}_{i+1}^1)$ is 
annihilated by the time one gets to $x_0 \in \s{C}_{i+1}$.  

Conversely, suppose that one starts with an element $x \in \ker
\Phi_{i+1}(0)/$ $\im \Phi_{i+2}(0)$.  Choose a 
representative $x_0 \in
\ker \Phi_{i+1}(0)$ for $x$.  Let $\eta_0$ be a preimage of $x_0$ in
$\s{F}_{i+1}^0$.  By the commutativity of \eqref{55}, there exists
$\zeta_0 \in \s{F}_i^1$ such that $\partial_1 \zeta_0 = \Phi_{i+1}
\eta_0$.  Hence this $\zeta_0$ is part of a homology class $\zeta \in
\ker \partial'_1/\im \partial'_2$.  We claim that $\zeta$ depends only
on the choice of $x$.  Indeed, by the exactness of the middle row in
\eqref{55}, a different choice of $\eta_0 \in \s{F}_{i+1}^0$
corresponds to a perturbation by an element in $\partial_1
(\s{F}_{i+1}^1)$.  By the commutativity of \eqref{55}, this 
corresponds to a perturbation of $\zeta_0$ by an element of 
$\Phi_{i+1} (\s{F}_{i+1}^1)$, yielding the same homology class 
$\zeta$.  A different choice of $\zeta_0$ so that $\partial_1 \zeta_0 
= \Phi_{i+1}
\eta_0$ corresponds, by exactness of the bottom row of \eqref{55}, to 
a perturbation by an element in $\partial_{2} (\s{F}_i^2)$ which 
obviously yields the same homology class $\eta$.  Finally, a different 
choice of $x_0$ corresponds to a perturbation by an element in 
$\Phi_{i+2}(0)(\s{C}_{i+2})$.  By the exactness of the top row and 
second to last column of \eqref{55}, one sees that this perturbation 
is annihilated by the time one arrives at $\s{F}_i^0$.  

Clearly the above two maps are inverses of each other.  Hence the 
desired isomorphism is established.  

\end{proof}

We can now establish the statement of the theorem.  Labeling more 
explicitly now, we let 

\begin{equation*} 
\begin{CD}
\cdots @>\partial_2^i>> \s{H}_i^1 @>\partial_1^i>> \s{H}_i^0 @>>> 0
\end{CD}
\end{equation*}

\noindent be the Koszul complex of $\s{H}_i$.  For $k 
\geq 2$, Claim~\ref{51} yields 
the following sequence of isomorphisms:

\begin{equation} \label{57}
\frac{\ker \partial_k^0}{\im \partial_{k+1}^0} \cong \frac{\ker 
\partial_{k-1}^{1}}{\im \partial_k^{1}} \cong \cdots \frac{\ker 
\partial_1^{k-1}}{\im \partial_2^{k-1}}.
\end{equation}

\noindent We then use Claim~\ref{54} to finish off the sequence:

\begin{equation} \label{58}
\frac{\ker \partial_1^{k-1}}{\im \partial_2^{k-1}} \cong 
\frac{\ker \Phi_{k}(0)}{\im \Phi_{k+1}(0)}.
\end{equation}

\noindent Together, \eqref{57} and \eqref{58} establish the 
statement of the 
theorem.  

\end{proof}

Our goal is to describe the homology of the localized complex 

\begin{equation*} 
\begin{CD}
\cdots @>\Phi_3(\lambda)>> \s{C}_2 @>\Phi_2(\lambda)>> \s{C}_1 
@>\Phi_1(\lambda)>> \s{C}_0.
\end{CD}
\end{equation*}

\noindent for an arbitrary $\lambda \in B_d$.  We will attain this 
goal in Section~\ref{85}, but we need some 
machinery first.  This is the subject of the next section.  

\section{The M\"obius Transform} \label{60}

In this section we define the notion of a M\"obius transform of a pure 
contractive $\s{A}_d$-module.  We begin by summarizing the main 
properties of M\"obius transforms on the unit ball in $\mbb{C}^d$.  
For a more detailed exposition, we refer the reader to 
\cite{rudinftotub}.  Recall that a M\"obius transform on the unit 
ball in $\mbb{C}^d$ is 
a 
continuous bijection 
$\varphi: \overline{B_d} \lra \overline{B_d}$ which satisfies the 
following (somewhat redundant) properties:

\begin{enumerate}
\item \label{60.5}
$\varphi$ is holomorphic in $B_d$.
\item \label{60.75}
$\varphi(B_d) = B_d$.
\item \label{60.825}
$\varphi(\partial B_d) = \partial B_d$.
\end{enumerate}

\noindent Let $\lambda \in B_d$, and define $\varphi_\lambda$ as 
follows:

\begin{equation} \label{61}
\varphi_\lambda(z_1, z_2, \ldots, z_d) = 
\left(\frac{\sqrt{1-|\lambda|^2} 
z_1}{1-\sum_{k=1}^d \, \overline{\lambda_k} z_k} - \lambda_1, \ldots, 
\frac{\sqrt{1-|\lambda|^2}
z_d}{1-\sum_{k=1}^d \, \overline{\lambda_k} z_k} - \lambda_d \right)
\end{equation}

\noindent or in ``vector notation'', 

\begin{equation*} 
\varphi_\lambda (z) = \frac{P_{\mbb{C}\lambda} (z - \lambda) - 
(1-|\lambda|^2) 
P_{\mbb{C} \lambda}^\perp z}{1-\left< z, \lambda\right>}, \; z \in 
B_d, 
\end{equation*}

\noindent where $P_{\mbb{C}\lambda}$ is the orthogonal projection onto 
the 
one-dimensional subspace $\mbb{C} \lambda \appa \mbb{C}^d$.  Some 
calculations reveal that $\varphi_\lambda$ is a M\"obius transform and 
that $\varphi_\lambda(\lambda) = 0$.  It is a striking fact that {\it 
any} 
M\"obius transform $\varphi$ can be written in the form $u \circ 
\varphi_\lambda$, where $u$ is a unitary operator on $\mbb{C}^d$ and 
$\lambda = \varphi^{-1}(0)$.  This fact allows us to find a useful 
formula for the expression $\left< \varphi(w), \varphi(z) \right>$, 
where $z, w \in B_d$.  Calculating this first for the case where 
$\varphi = \varphi_\lambda$, we obtain the following identity:

\begin{equation} \label{67}
\left< \varphi(z), \varphi(w) \right> = 1 
- 
\frac{(1-|\lambda|^2)(1-\left<z,w 
\right>)}{(1-\left<\lambda, w 
\right>)(1-\left< z, \lambda \right>)}, \; z, 
w \in B_d.
\end{equation}

\noindent Consequently, for any M\"obius transform $\varphi$ with 
$\lambda = \varphi^{-1}(0)$, $\left< \varphi(w), \varphi(z) 
\right>$ is also given by \eqref{67}.  

The following theorem unveils the role that is played by M\"obius 
transforms 
in 
the theory of free contractive $\s{A}_d$-modules.  

\begin{thm} \label{68}
Let $\varphi^1$, $\varphi^2$, $\ldots$, 
$\varphi^d$ be the coordinates of the M\"obius transform 
$\varphi$.  Define the map $\Phi: \overbrace{H_d^2 
\oplus 
\cdots H_d^2}^{d \text{~times~}} \lra H_d^2$ to be left 
multiplication 
by the row vector 
$(\varphi^1 \; \cdots \; \varphi^d)$, i.e. 

\begin{equation*} 
\Phi \left( \begin{array}{c}
\xi_1 \\
\xi_2 \\
\vdots \\
\xi_d \end{array} \right) = \sum_{k=1}^d \, \varphi^k \xi_k, 
\; 
\xi_1, 
\xi_2, \ldots, \xi_d \in H_d^2.
\end{equation*}

\noindent Then $\Phi$ is a partially isometric module 
homomorphism 
with range $\{\xi \in H_d^2: \xi(\lambda) = 0\} = 
\{k_\lambda\}^\perp$.  
\end{thm}

\begin{proof}
It is obvious that $\Phi$ is a module homomorphism.  To show 
that it 
is partially isometric with the stated range, it suffices to show that

\begin{equation*} \label{70}
\left<(1-\Phi \Phi^*) k_w, k_z \right> = \left< 
P_\lambda k_w, k_z 
\right>,
\end{equation*}

\noindent for any $w, z \in B_d$.  The sufficiency of this condition 
follows from the fact that the set of all $k_z$'s forms a spanning set 
of $H_d^2$.  

We compute the left side of \eqref{70}.  Using 
\eqref{3}, \eqref{8}, and formula \eqref{67}, we have 

\begin{multline} \label{71}
\left< (1-\Phi \Phi^*) k_w, k_z \right> = \frac{1}{1-\left<z, w 
\right>} - \frac{\left<\varphi(z), \varphi(w) 
\right>}{1 - \left< z, w \right>} \\
= \frac{1}{1-\left<z, w
\right>} - \frac{1}{1-\left<z, w
\right>} + \frac{1-|\lambda|^2}{(1 - \left< \lambda, w \right>)(1 - 
\left< z, \lambda \right>)} \\
= \frac{1-|\lambda|^2}{(1 - \left< 
\lambda, w \right>)(1 -
\left< z, \lambda \right>)}.
\end{multline}

\noindent We compute the right side of \eqref{70} as follows:

\begin{equation} \label{72}
\left< P_\lambda k_w, k_z \right> = \frac{\left<k_w, k_\lambda 
\right>}{\|k_\lambda\|^2} \left< k_\lambda, k_z \right> 
= \frac{1-|\lambda|^2}{(1-\left< \lambda, w \right>)(1-\left<z, 
\lambda \right>)},
\end{equation}

\noindent which is identical to \eqref{71}, hence we have established 
\eqref{70}.  

\end{proof}

We are now in a position to define a M\"obius transform of $H_d^2$.  

\begin{mydef} \label{73}
Let $\varphi$ be a M\"obius transform.  We define $(H_d^2)_\varphi$ 
to be 
the $\s{A}_d$-module whose underlying Hilbert space is $H_d^2$ and 
whose $\s{A}_d$ is given by $z_i \cdot \xi = \varphi^i \xi$ for any 
$\xi \in H_d^2$ and $i = 1, 2, \ldots, d$, where $\varphi^i$ is the 
$i$th coordinate function of $\varphi$.  
\end{mydef}

\begin{thm} 
Let $\varphi$ be a M\"obius transform.  Then $(H_d^2)_\varphi$ is a 
pure contractive $\s{A}_d$-module.  Furthermore, if we set $\lambda = 
\varphi^{-1}(0)$, then the map $U_{\varphi}: H_d^2 
\lra 
(H_d^2)_\varphi$ defined by $U_{\varphi} \xi = (\xi \circ 
\varphi) 
\frac{k_\lambda}{\|k_\lambda\|}$ is a unitary module isomorphism.  
\end{thm}

\begin{proof}
Let $(\Phi_1, \Phi_2, \ldots, \Phi_d)$ be the 
$d$-tuple such that $\Phi_i \in \s{B}(H_d^2)$ is multiplication by 
$\varphi^i$.  By conditions \eqref{60.5}, \eqref{60.75}, and 
\eqref{60.825}, to prove the first part of the theorem it suffices 
to show that the 
following two 
conditions on 
$(\Phi_1, 
\Phi_2, \ldots, \Phi_d)$ hold.
\begin{enumerate}
\item \label{75}
$\sum_{k=1}^d \, \Phi_k \Phi_k^* \leq 1$.
\item \label{76}
$\wotlim_{n \sra \infty} \, \sum_{i_1, i_2, \ldots, i_n}^d \, 
\Phi_{i_1} \Phi_{i_2} \cdots \Phi_{i_n} 
\Phi_{i_n}^* \cdots 
\Phi_{i_2}^* 
\Phi_{i_1}^* = 0$.  
\end{enumerate}

\noindent The first condition \eqref{75} follows from 
Theorem~\ref{68}, since the row operator $(\Phi_1, 
\Phi_2, \cdots, \Phi_d)$ is a partial isometry.  For the second 
condition \eqref{76}, we apply the stated limit to the linear 
functional $\left< 
(\cdot) k_z, k_w \right>$:

\begin{equation} \label{77}
\lim_{n \sra \infty} \sum_{i_1, i_2, \ldots, i_n}^d \, \left< 
\Phi_{i_1} \cdots \Phi_{i_n} \Phi_{i_n}^* 
\cdots \Phi_{i_1}^* k_z, k_w \right> = \lim_{n \sra \infty} 
\, \left< \varphi_\lambda(w), \varphi_\lambda(z) \right>^n \left< k_z, 
k_w \right>.
\end{equation}

\noindent where we make use of \eqref{8}.  By property 
\eqref{60.75} of the classical M\"obius transform, 
$|\left< \varphi_\lambda(w), \varphi_\lambda(z) \right>| < 1$.  Hence 
the limit in \eqref{77} tends to $0$.  Since the set of all $k_z$ form 
a spanning set for $H_d^2$, and since the sums $\sum_{i_1, i_2, 
\ldots, i_n}^d \, \Phi_{i_1} \cdots \Phi_{i_n} 
\Phi_{i_n}^* \cdots \Phi_{i_1}^*$ are uniformly 
bounded over $n$, it follows that condition \eqref{76} is valid.  

For the next statement in the proof, recall that by 
Theorem~\ref{5.5} 
there exists a minimal dilation $U_\varphi: H_d^2 \otimes \s{D} \lra 
(H_d^2)_\varphi$.  Let $(S'_1, S'_2, \ldots, S'_d)$ be the $d$-tuple 
associated with $H_d^2 \otimes \s{D}$.  By Theorem~\ref{68} 
$\sum_{k=1}^d \, S'_k S_k^{' \; *}$ is the projection onto $(k_0 
\otimes \s{D})^\perp = \{\xi \in H_d^2 \otimes \s{D} : \xi(0) = 0\}$.  
Since $U_\varphi$ is a coisometric module homomorphism, we have 
the 
following 
equation:  Let $P'_0$ be the orthogonal projection onto the 
space 
$k_0 \otimes \s{D}$.  

\begin{equation} \label{78}
U_\varphi P'_0 U_\lambda^* = U_\varphi (1 -\sum_{k=1}^d \, S'_k 
S_k^{' 
\; *}) U_\varphi^* = 
1-\sum_{k=1}^d \, \Phi_i \Phi_i^* = 
P_\lambda
\end{equation}

\noindent Theorem~\ref{68}.  Hence $P'_0 U_\varphi^*$ is 
a rank one operator, and by the 
minimality 
condition on $U_\varphi$, the submodule generated by the image 
this 
operator 
must be $H_d^2 \otimes \s{D}$.  It follows that $\dim \s{D} = 1$, 
hence we may view $U_\varphi$ as a map $H_d^2 \lra 
(H_d^2)_\varphi$.  
The 
equation in \eqref{78} also implies that $U_\varphi P_0$ is a partial 
isometry with support $\mbb{C} k_0$ and image $\mbb{C} k_\lambda$.  
Hence we may set $U_\varphi 1 = \frac{k_\lambda}{\|k_\lambda\|}$.  The 
equation 
$U_\varphi \xi = 
(\xi \circ \varphi)\frac{k_\lambda}{\|k_\lambda\|}$ now 
follows since $U_\varphi$ is a module homomorphism.  

To complete the proof, it suffices to show that $\ker U_\varphi 
= 
\{0\}$.  
To this end, let $\xi \in H_d^2$, and suppose that $(\xi \circ 
\varphi) \frac{k_\lambda}{\|k_\lambda\|} = 0$.  Since 
$k_\lambda(z)/\|k_\lambda\| = \frac{\sqrt{1-|\lambda|^2}}{1-\left< z, 
\lambda \right>}$ is never $0$ on $B_d$, we must have $\xi 
(\varphi(z))$ for all $z \in B_d$.  Hence $\xi(z) = 0$ for all 
$z \in B_d$ since $\varphi$ is bijective on $B_d$.  Hence $\xi = 
0$.  
\end{proof}

The following corollary shows how a M\"obius transform 
$U_\varphi$ provides a means of ``changing the base point'' 
from $0$ to $\lambda \in B_d$ when considering the module 
$H_d^2$.  

\begin{cor} 
Let $\varphi$ be a M\"obius transform and let $\lambda = 
\varphi^{-1}(0)$.  Then $U_{\varphi} \{k_0\}^\perp = 
\{k_\lambda\}^\perp$.  
\end{cor}

\begin{proof}
This is obvious from the definition of $U_{\varphi}$.  

\end{proof}

To conclude this section, we demonstrate an ``ergodicity'' 
property of the set of M\"obius transforms.  

\begin{thm} 
Let $\s{M}$ be a proper non-trivial closed submodule of $H_d^2$.  
Then there exists $\xi \in \s{M}$ and a M\"obius transform $\varphi$
such that $U_\varphi \xi \notin \s{M}$.  
\end{thm}

\begin{proof}
Since $\s{M}$ is a proper closed submodule, it cannot contain $k_0$.  
Hence 

\begin{equation*} 
M = \sup 
\{ |\left< \xi, k_0 \right> | : \|\xi\| = 1, \; \xi \in \s{M} \} < 1.
\end{equation*}

\noindent By Theorem~3.2 in \cite{g1}, 
there 
exists $\lambda \in B_d$ such that 

\begin{equation} \label{82}
\frac{\|P_\s{M} k_\lambda\|}{\|k_\lambda\|} > M.
\end{equation}

\noindent An explicit calculation involving \eqref{61} shows that 

\begin{equation*} 
U_{\varphi_{-\lambda}} U_{\varphi_\lambda} \xi = \xi \circ u,
\end{equation*}

\noindent where $u$ is a unitary operator on $\mbb{C}^d$.  Hence by 
the definition of $U_u$, $U_{\varphi_{-\lambda}} U_{\varphi_\lambda} = 
U_u$.  Hence $U_{\varphi_\lambda}^* = U_{u^*} U_{\varphi_{-\lambda}}$.  
Therefore by our assumption that $\s{M}$ is invariant under M\"obius 
transforms, it follows that 

\begin{equation} \label{84}
\sup \{ |\left< U_{\varphi_\lambda}^* \xi, k_0 \right>| : \|\xi\| = 1, 
\; \xi \in \s{M} \} \leq M.
\end{equation}

\noindent By definition, $U_{\varphi_\lambda} k_0 = 
\frac{k_\lambda}{\|k_\lambda\|}$.  Hence by \eqref{84}, $|\left< \xi, 
\frac{k_\lambda}{\|k_\lambda\|} \right>| \leq M$ for all $\xi \in 
\s{M}$ such that $\|\xi\| = 1$.  But this implies that 

\begin{equation*} 
\frac{\|P_\s{M} k_\lambda\|}{\|k_\lambda\|} \leq M, 
\end{equation*}

\noindent which contradicts \eqref{82}.  

\end{proof}

\section{The homology of localized free resolutions} \label{85}

In this section we use the machinery developed in 
Section~\ref{60} to
extend Theorem~\ref{39.25}.  We first generalize 
Definition~\ref{73}.

\begin{mydef} 
Let $\s{H}$ be a pure contractive $\s{A}_d$-module, and let $\varphi$ 
be a M\"obius transform.  We define the M\"obius transform of $\s{H}$ 
by $\varphi$ to be the module $(\s{H})_\varphi$ with underlying 
Hilbert space $\s{H}$ and $\s{A}_d$ action defined by $z_i \cdot \xi = 
\varphi^i(T_1, T_2, \ldots, T_d) \xi$ for all $\xi \in H_d^2$ and $i = 
1, 2, \ldots, d$.  
\end{mydef}

\begin{thm} 
Let $\s{H}$ be a pure contractive $\s{A}_d$-module, and let $\varphi$ 
be a M\"obius transform.  Then $(\s{H})_\varphi$ is contractive and 
pure.  
\end{thm}

\begin{proof}
Let $V: H_d^2 \otimes \s{D} \lra \s{H}$ be a dilation of $\s{H}$.  
From Definition~\ref{73}, we see that $V$ is also a module 
homomorphism as a map from $(H_d^2)_\varphi \otimes \s{D}$ onto 
$(\s{H})_\varphi$.  Precomposing $V$ with $U_\varphi \otimes I_\s{D}$ 
gives a coisometric module homomorphism $V': H_d^2 \otimes \s{D} \lra 
(\s{H})_\varphi$.  Hence $(\s{H})_\varphi$ is isomorphic to the 
quotient of a free module, whence it is pure and contractive.  
\end{proof}

The main result of this section is the following:

\begin{thm} \label{87}
Let $\s{H}$ be a pure contractive $\s{A}_d$-module, and let 
$\varphi$ be a M\"obius transform with $\lambda = \varphi^{-1}(0)$.  
Let 

\begin{equation*} \label{88}
\begin{CD}
\cdots @>\partial'_3>> E_2 @>\partial'_2>> E_1 @>\partial'_1>> E_0 
@>>> 0
\end{CD}
\end{equation*}

\noindent be the Koszul complex of $(\s{H})_\varphi$, and let 

\begin{equation} \label{89}
\begin{CD}
\cdots @>\Phi_3(\lambda)>> \s{C}_2 @>\Phi_2(\lambda)>> \s{C}_1 
@>\Phi_1(\lambda)>> \s{C}_0
\end{CD}
\end{equation}

\noindent be the localization at $\lambda$ of a free 
resolution of 
$\s{H}$.  Then for $k \geq 1$,

\begin{equation*} 
\frac{\ker \partial'_k}{\im \partial'_{k+1}} \cong \frac{\ker 
\Phi_k(\lambda)}{\im \Phi_{k+1}(\lambda)}.
\end{equation*}

\end{thm}

\begin{proof}
Let 

\begin{equation} \label{91}
\begin{CD}
\cdots @>\Phi_2>> H_d^2 \otimes \s{C}_1 @>\Phi_1>> H_d^2 \otimes 
\s{C}_0 @>\Phi_0>> \s{H} @>>> 0
\end{CD}
\end{equation}

\noindent be the free resolution of $\s{H}$ from which \eqref{89} is
derived.  Since $\Phi_i$ is a module homomorphism it follows that for
$i \geq 0$, $\varphi^j \Phi_i = \Phi_i \varphi^j$ for $j = 1, 2,
\ldots, d$.  Hence we may view \eqref{91} as an exact sequence of 
M\"obius transformed modules:

\begin{equation*} 
\begin{CD}
\cdots @>\Phi_2>> (H_d^2)_\varphi \otimes \s{C}_1 @>\Phi_1>> 
(H_d^2)_\varphi \otimes \s{C}_0 @>\Phi_0>> (\s{H})_\varphi @>>> 0.
\end{CD}
\end{equation*}

\noindent For $i \geq 1$ let $\Phi'_i = U_\varphi^* \Phi_i U_\varphi$.  
Then we have the following partial isomorphism of complexes:

\begin{equation*} 
\begin{CD}
\cdots @>\Phi_2>> (H_d^2)_\varphi \otimes \s{C}_1 @>\Phi_1>>
(H_d^2)_\varphi \otimes \s{C}_0 @>\Phi_0>> (\s{H})_\varphi @>>> 0 \\
& &  @VU_\varphi^*VV @VU_\varphi^*VV \\
\cdots @>\Phi'_2>> H_d^2 \otimes \s{C}_1 @>\Phi'_1>> H_d^2 \otimes 
\s{C}_0.
\end{CD}
\end{equation*}

\noindent Since two isomorphic complexes have isomorphic homologies, 
there is a coisometric homomorphism $\Phi'_0: H_d^2 \otimes \s{C}_0 
\lra (\s{H})_\varphi$ which make the following sequence exact:

\begin{equation*} 
\begin{CD}
\cdots @>\Phi'_2>> H_d^2 \otimes \s{C}_1 @>\Phi'_1>> H_d^2 \otimes 
\s{C}_0 @>\Phi'_0>> (\s{H})_\varphi @>>> 0.
\end{CD}
\end{equation*}

\noindent This is obviously a free resolution for $(\s{H})_\varphi$, 
hence Theorem~\ref{39.25} implies that 

\begin{equation*}
\frac{\ker \partial'_k}{\im \partial'_{k+1}} \cong \frac{\ker 
\Phi'_k(0)}{\im \Phi'_{k+1}(0)}
\end{equation*}

\noindent for $k \geq 1$.  We now compute $\Phi'_i(0)$ for $i \geq 
1$.  Let $\eta 
\in \s{C}_i$ and $\eta' \in \s{C}_{i+1}$.  Then 

\begin{multline*} 
\left< \Phi'_i(0) \eta, \eta' \right> = \left< \Phi'_i (k_0 \otimes 
\eta), k_0 \otimes \eta' \right> = \left< U_\varphi^* \Phi_i U_\varphi 
(k_0 \otimes \eta), k_0 \otimes \eta' \right> \\
= \left<\Phi_i (\frac{k_\lambda}{\|k_\lambda\|} \otimes \eta), 
\frac{k_\lambda}{\|k_\lambda\|} \otimes \eta' \right> = 
\left< \frac{k_\lambda}{\|k_\lambda\|} \otimes \eta, \Phi_i^* 
(\frac{k_\lambda}{\|k_\lambda\|} \otimes \eta') \right> \\
= \left< \frac{k_\lambda}{\|k_\lambda\|} \otimes \eta, 
\frac{k_\lambda}{\|k_\lambda\|} \otimes \Phi_i (\lambda)^* \eta' 
\right> 
= \frac{\|k_\lambda\|^2 \left< \eta, \Phi_i(\lambda)^* \eta' 
\right>}{\|k_\lambda\|^2} = \left< \Phi_i(\lambda) \eta, \eta' 
\right>.
\end{multline*}

\noindent Since $\eta$ and $\eta'$ were arbitrary, $\Phi'_i(0) = 
\Phi_i(\lambda)$ for each $i \geq 1$.  
The conclusion of the theorem now follows.  
\end{proof}



\providecommand{\bysame}{\leavevmode\hbox to3em{\hrulefill}\thinspace}
\providecommand{\MR}{\relax\ifhmode\unskip\space\fi MR }
\providecommand{\MRhref}[2]{%
  \href{http://www.ams.org/mathscinet-getitem?mr=#1}{#2}
}
\providecommand{\href}[2]{#2}


\end{document}